\documentclass[graybox]{svmult}%

\usepackage{mathptmx}%
\usepackage{helvet}%
\usepackage{courier}%
\usepackage{makeidx}%
\usepackage{graphicx}%
\usepackage{multicol}%
\usepackage{footmisc}%
\usepackage{subfigure}%
\usepackage{amsmath,amssymb}%

\newcommand{\N}{\mathbb{N}}%
\newcommand{\Z}{\mathbb{Z}}%
\newcommand{\R}{\mathbb{R}}%
\newcommand{\tp}{\mathrm{top}}%
\newcommand{\ep}{\varepsilon}%
\newcommand{\AC}{\mathcal{A}}%
\newcommand{\LC}{\mathcal{L}}%
\newcommand{\UC}{\mathcal{U}}%
\newcommand{\PC}{\mathcal{P}}%
\newcommand{\RC}{\mathcal{R}}%
\newcommand{\QC}{\mathcal{Q}}%
\newcommand{\FC}{\mathcal{F}}%
\newcommand{\EC}{\mathcal{E}}%
\newcommand{\id}{\mathrm{id}}%
\newcommand{\rmd}{\mathrm{d}}%
\newcommand{\rmD}{\mathrm{D}}%
\newcommand{\rmS}{\mathrm{S}}%
\newcommand{\Mis}{\mathrm{M}}%
\newcommand{\diam}{\mathrm{diam}}%
\newcommand{\vol}{\mathrm{vol}}%

\begin{document}

\title*{Entropy of nonautonomous dynamical systems}%
\author{Christoph Kawan}%
\institute{Fakult\"{a}t f\"{u}r Informatik und Mathematik, Universit\"{a}t Passau, Innstra{\ss}e 33, 94032 Passau, Germany; e-mail: christoph.kawan@uni-passau.de}%
\maketitle%

\abstract{Different notions of entropy play a fundamental role in the classical theory of dynamical systems. Unlike many other concepts used to analyze autonomous dynamics, both measure-theoretic and topological entropy can be extended quite naturally to discrete-time nonautonomous dynamical systems given in the process formulation. This paper provides an overview of the author's work on this subject. Also an example is presented that has not appeared before in the literature.%
}%

\section{Introduction}%

In the 1950s, Kolmogorov and Sinai established the concept of measure-theoretic (or metric) entropy, based on Shannon entropy from information theory, as an invariant for measure-preserving maps on probability spaces. This invariant was used, e.g., by Ornstein \cite{Orn} to classify Bernoulli shifts. Some years later, Adler, Konheim and McAndrew \cite{AKM} defined in strict analogy a notion of entropy for continuous maps on compact spaces. They already conjectured that both entropy notions are related to each other in the sense of a variational principle, i.e., the topological entropy equals the supremum over all measure-theoretic entropies (supremizing over all invariant Borel probability measures). This was proved not much later by Goodman, Goodwyn and Dinaburg \cite{Go1,Gow,Din}.%

In the theory of dynamical systems, developed in the ensuing decades, both notions of entropy play a fundamental role as it turned out that they are related to many other dynamical characteristics such as Lyapunov exponents, dimensions of invariant measures and invariant sets and growth rates of periodic orbits, but also to the existence of horseshoes. Moreover, entropy has become a central concept in a branch of the topological theory of dynamical systems dedicated to the question of how well a dynamical system can be `digitalized', i.e., modeled by a symbolic dynamical system \cite{Dow}.%

Motivated by the study of triangular maps, Kolyada and Snoha \cite{KSn} extended the notion of topological entropy to nonautonomous systems given by a sequence of continuous maps on a compact metric space. Together with Misiurewicz, they generalized this concept to sequences of maps between possibly different metric spaces in \cite{KMS} and proved analogues of the Misiurewicz-Szlenk formula for the entropy of piecewise monotone interval maps. Further work on topological entropy of nonautonomous systems has been done in \cite{PMa,Mou,OW1,OW2,ZCh,ZZH,ZLX} by several researchers with different motivations and partially independently of \cite{KSn,KMS}. An essential difference to the classical theory that should be mentioned is that the nonautonomous version of topological entropy is \emph{not} a purely topological quantity. In fact, it depends on the sequence of metrics imposed on the time-varying state space.%

Concepts of measure-theoretic entropy for sequences of maps were first introduced in the papers \cite{ZLX,Can,Ka1}. While \cite{ZLX,Can} require that all maps in the sequence preserve the same measure, a very restrictive condition, the approach in \cite{Ka1} is completely general. The invariant measure now becomes a sequence $(\mu_n)_{n\in\Z_+}$ of measures so that $(f_n)_*\mu_n = \mu_{n+1}$ for the given sequence of maps $f_n$. To introduce a reasonable notion of entropy in this general context, an additional structure (called an \emph{admissible class}) needs to be imposed on the system, consisting in a family of sequences of measurable partitions. This family has to satisfy certain axioms in order to obtain structural results such as a power rule and invariance under a reasonably general class of transformations.%

In the topological framework, a relation between the topological and the measure-theoretic entropy can be established through the definition of a suitable admissible class adapted to the metric space structure. We call this class the \emph{Misiurewicz class}, since it allows for an easy adaptation of Misiurewicz's proof of the variational principle \cite{Mis} to show that the measure-theoretic entropy is bounded above by the topological entropy. In the classical case of a single map, the entropy computed with respect to the Misiurewicz class reduces again to the Kolmogorov-Sinai measure-theoretic entropy.%

It is still unclear whether a full variational principle holds in this context. One obstruction to a proof, amongst others, is that the Misiurewicz class might not contain elements of arbitrarily small diameter, in general. Some sufficient conditions for the existence of such sequences of small-diameter partitions have been identified in \cite{KLa}, but a general approach to this problem is still missing.%

The paper is organized as follows. In Section \ref{sec_motiv}, we motivate the entropy theory for nonautonomous dynamical systems by applications in networked control. Section \ref{sec_entth} explains the entropy theory developed in \cite{KSn,KMS} and \cite{Ka1,Ka2,KLa}, including the nonautonomous versions of topological and measure-theoretic entropy and their relation. Finally, an example for a system satisfying a full variational principle is presented in Section \ref{sec_ex}.%

\section{Motivation from networked control}\label{sec_motiv}

\begin{figure}[b]
\begin{center}
\includegraphics[scale=.65]{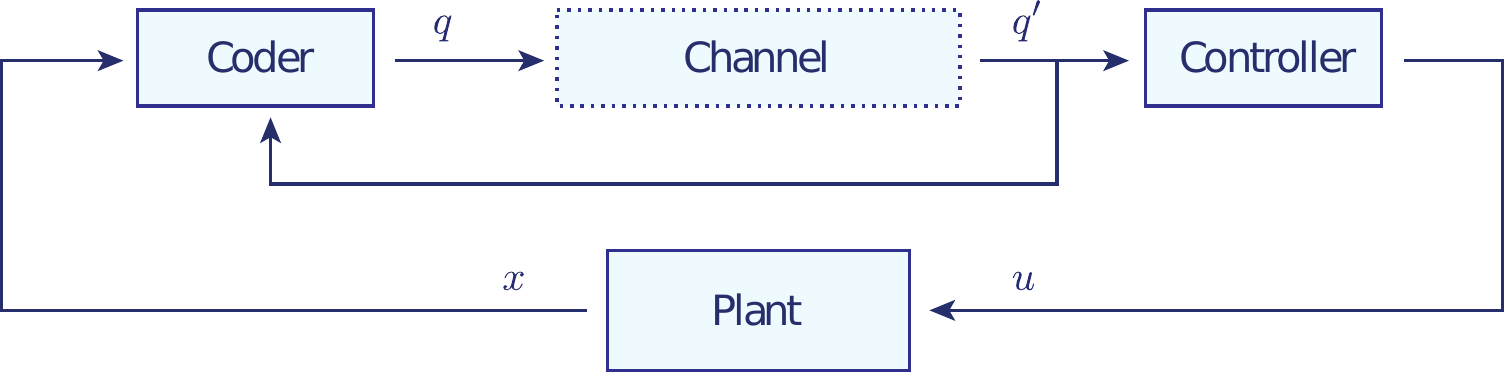}
\end{center}
\caption{The simplest model of an NCS}
\label{fig:0}
\end{figure}

The author's central motivation for the development of a nonautonomous entropy theory comes from problems arising in networked control. Networked control systems (NCS) are spatially distributed systems whose components (sensors, controllers and actuators) share a common digital communication network. Examples can be found in vehicle tracking, underwater communications for remotely controlled surveillance and rescue submarines, remote surgery, space exploration and aircraft design. Another large field of applications can be found in modern industrial systems, where industrial production is combined with information and communication technology (`Industry 4.0'). A fundamental problem in this field is to determine the minimal requirements on the communication network so that a specified control objective can be achieved.%

The simplest model of an NCS consists of a single feedback loop containing a finite-capacity channel which transmits state information acquired by a sensor from a coder to the controller (see Fig.~\ref{fig:0}). The first task of the controller, before deciding on the control action, often consists in the computation of a state estimate. If the system is autonomous, it has been shown in \cite{MPo} that the smallest channel capacity above which a state estimation of arbitrary precision can be achieved is given by the topological entropy of the system. If the problem setup is slightly changed, time-dependencies of many different sorts can appear. Here are some examples:%
\begin{itemize}
\item Non-invariance of the region of relevant initial states leads to a time-dependent state space.%
\item The requirement of an exponential improvement of the estimate over time leads to a time-dependent metric on the state space.%
\item In a stochastic formulation of the problem, non-invariance of the distribution of $x_0$ (the initial state) leads to a time-dependent probability measure.%
\item Time-varying coding policies lead to time-dependent partitions of the state space (with respect to which entropy needs to be computed).%
\end{itemize}
The entropy theory described in this paper is sufficiently general to handle all of these time-dependencies. A first application to a state estimation problem can be found in \cite{Ka3}.%

\section{Entropy theory for nonautonomous systems}\label{sec_entth}

{\bf Notation:} We write $\N = \{1,2,3,\ldots\}$ and $\Z_+ = \{0,1,2,\ldots\}$. By $\delta_x$ we denote the Dirac measure concentrated at a point $x$. The cardinality of a finite set $S$ is denoted by $\# S$. If $A$ is a subset of a metric space $(X,d)$, we write $\diam A = \sup \{ d(x,y) : x,y\in A\}$. If $\AC$ is a collection of sets $A\subset X$, we write $\diam\AC = \sup\{\diam A: A\in\AC\}$. All logarithms are taken to the base $2$.%

A nonautonomous dynamical system, or briefly an NDS, is a pair $(X_{\infty},f_{\infty})$, where $X_{\infty} = (X_n)_{n\in\Z_+}$ is a sequence of sets and $f_{\infty} = (f_n)_{n\in\Z_+}$ a sequence of maps $f_n:X_n \rightarrow X_{n+1}$. For all $i\in\Z_+$ and $n\in\N$, we define%
\begin{equation*}
  f_i^0 := \id_{X_i},\quad f_i^n := f_{i+n-1} \circ \cdots \circ f_{i+1} \circ f_i,\quad f_i^{-n} := (f_i^n)^{-1}.%
\end{equation*}
We do not assume that the maps $f_i$ are invertible, so $f_i^{-n}$ is only applied to sets. We speak of a \emph{topological NDS} if each $X_n$ is a compact metric space $(X_n,d_n)$ and the sequence $f_{\infty}$ is equicontinuous, i.e., for every $\ep>0$ there is $\delta>0$ so that $d_n(x,y) < \delta$ for any $n\in\Z_+$ and $x,y\in X_n$ implies $d_{n+1}(f_n(x),f_n(y)) < \ep$.%

\subsection{Topological entropy}\label{subsec_topent}

To define the topological entropy of a dynamical system, one needs to specify a \emph{resolution} on the state space. Usually, this resolution is given by a finite $\ep>0$ or by an open cover. In the case of an NDS $(X_{\infty},f_{\infty})$, we have to consider a sequence of open covers instead. Hence, let $\UC_{\infty} = (\UC_n)_{n\in\Z_+}$ be a sequence so that $\UC_n$ is an open cover of $X_n$ for every $n$. For all $i\in\Z_+$ and $n\in\N$ define%
\begin{equation*}
  \UC_i^n := \bigvee_{j=0}^{n-1}f_i^{-j}\UC_{i+j},%
\end{equation*}
which is the common refinement of the open covers $f_i^{-j}\UC_{i+j}$ of $X_i$, i.e., the open cover whose elements are of the form%
\begin{equation*}
  U_{j_i} \cap f_i^{-1}(U_{j_{i+1}}) \cap \ldots \cap f_i^{-n+1}(U_{j_{i+n-1}}),\quad U_{j_l} \in \UC_l.%
\end{equation*}
Then the entropy of $f_{\infty}$ w.r.t.~$\UC_{\infty}$ is defined by%
\begin{equation}\label{eq_defent_ocs}
  h(f_{\infty};\UC_{\infty}) := \limsup_{n\rightarrow\infty}\frac{1}{n}\log N(\UC_0^n),%
\end{equation}
where $N(\cdot)$ denotes the minimal cardinality of a finite subcover. Here, unlike in the autonomous case, the $\limsup$ in general is not a limit (see \cite{KSn} for a counter-example).%

To define a notion of topological entropy, independent of a given resolution, one usually takes the supremum over all resolutions. However, taking the supremum of $h(f_{\infty};\UC_{\infty})$ over all sequences $\UC_{\infty}$ would result in a quantity that is usually $+\infty$, because a sequence of open covers whose diameters exponentially converge to zero generates an increase of information that is not due to the dynamics of the system. Hence, such sequences have to be excluded. An elegant way how to do this, is to consider only sequences with Lebesgue numbers bounded away from zero. We thus let $\LC(X_{\infty})$ denote the family of all such sequences and define the \emph{topological entropy of $(X_{\infty},f_{\infty})$} as%
\begin{equation*}
  h_{\tp}(f_{\infty}) := \sup_{\UC_{\infty}\in\LC(X_{\infty})}h(f_{\infty};\UC_{\infty}).%
\end{equation*}
This definition was first given in \cite{KMS}. Some properties of $h_{\tp}$ are the following:%
\begin{itemize}
\item Alternative characterizations in terms of $(n,\ep)$-spanning or $(n,\ep)$-separated sets can be given. For instance, a set $E \subset X_0$ is $(n,\ep;f_{\infty})$-spanning if for every $x\in X_0$ there exists $y\in E$ such that $d_i(f_0^i(x),f_0^i(y)) < \ep$ for $0 \leq i < n$. Letting $r(n,\ep;f_{\infty})$ denote the minimal cardinality of an $(n,\ep;f_{\infty})$-spanning set,%
\begin{equation}\label{eq_char_spansets}
  h_{\tp}(f_{\infty}) = \lim_{\ep\downarrow0}\limsup_{n\rightarrow\infty}\frac{1}{n}\log r(n,\ep;f_{\infty}).%
\end{equation}
\item In the case where $X_{\infty}$, $d_{\infty}$ and $f_{\infty}$ are constant, $h_{\tp}(f_{\infty})$ reduces to the usual notion of topological entropy for maps, which immediately follows from \eqref{eq_char_spansets}.%
\item The topological entropy $h_{\tp}(f_{\infty})$ also generalizes several other notions of entropy studied before, as for instance \emph{topological sequence entropy} \cite{Go2} and \emph{topological entropy for uniformly continuous maps on non-compact metric spaces} \cite{Bow}.%
\item Fundamental properties of topological entropy for maps carry over to its nonautonomous generalization, as for instance the power rule, which can be formulated as follows. For $m\in\N$ define the $m$th power system $(X_{\infty}^{[m]},f_{\infty}^{[m]})$ by $X_n^{[m]} := X_{nm}$ and $f_n^{[m]} := f_{nm}^m$. Then the following power rule holds:%
\begin{equation*}
  h_{\tp}(f_{\infty}^{[m]}) = m \cdot h_{\tp}(f_{\infty}).%
\end{equation*}
Here the equicontinuity of $f_{\infty}$ is essential, see \cite{KSn} for a counter-example in the case when $f_{\infty}$ is not equicontinuous.%
\end{itemize} 

\subsection{Measure-theoretic entropy}\label{subsec_mesent}%

To define measure-theoretic entropy, we consider systems given by measurable maps $f_n:X_n \rightarrow X_{n+1}$ between probability spaces $(X_n,\FC_n,\mu_n)$, preserving the measures $\mu_n$ in the sense that $(f_n)_*\mu_n = \mu_{n+1}$ for all $n\in\Z_+$. In this case, we also call the sequence $\mu_{\infty} = (\mu_n)_{n\in\Z_+}$ an \emph{invariant measure sequence}, or briefly an \emph{IMS} for the given NDS $(X_{\infty},f_{\infty})$, and we speak of a \emph{measure-theoretic NDS}. Analogously to the topological framework, we define the entropy of $f_{\infty}$ w.r.t.~a sequence of finite measurable partitions $\PC_n$ of $X_n$ by%
\begin{equation*}
  h(f_{\infty};\PC_{\infty}) = h_{\mu_{\infty}}(f_{\infty};\PC_{\infty}) := \limsup_{n\rightarrow\infty}\frac{1}{n}H_{\mu_0}(\PC_0^n),%
\end{equation*}
where $\PC_0^n$ denotes the partition $\bigvee_{i=0}^{n-1}f_0^{-i}\PC_i$ and $H_{\mu_0}(\cdot)$ is the Shannon entropy of a partition computed w.r.t.~the measure $\mu_0$.%

To define measure-theoretic entropy independently of a sequence of partitions, we have to follow a similar strategy as in the topological case. However, the concept of Lebesgue numbers is not helpful here, and a similar construction of a family $\LC(X_{\infty})$, using the measures $\mu_n$, does not lead to satisfying results. Looking at the topological theory, one sees that results for topological entropy such as the power rule rely on the equicontinuity of the sequence $f_{\infty}$, and not on the mere continuity of each $f_n$. However, in the measure-theoretic setting considered here we do not require a similar property.%

One way to overcome these obstructions is the study of the essential properties of the family $\LC(X_{\infty})$, defined in the topological framework, and enforcing these properties in the measure-theoretic framework by an axiomatic definition. As it turns out, the following definition leads to satisfying results.%

\begin{definition}
A nonempty family $\EC$ of sequences of finite measurable partitions for $X_{\infty}$ is called an \emph{admissible class} if it satisfies the following axioms:%
\begin{enumerate}
\item[(A)] For each $\PC_{\infty} = (\PC_n)_{n\in\Z_+}\in\EC$ there is a bound $N\in\N$ on the cardinality $\#\PC_n$, i.e., $\#\PC_n \leq N$ for all $n\in\Z_+$.%
\item[(B)] If $\PC_{\infty} = (\PC_n)_{n\in\Z_+} \in \EC$ and $\QC_{\infty} = (\QC_n)_{n\in\Z_+}$ is another sequence of finite measurable partitions for $X_{\infty}$ such that each $\QC_n$ is coarser than $\PC_n$, then $\QC_{\infty}\in\EC$.%
\item[(C)] If $\PC_{\infty} = (\PC_n)_{n\in\Z_+} \in \EC$ and $m\in\N$, then also the sequence $\PC_{\infty}^{\langle m \rangle}$, defined as follows, is an element of $\EC$:%
\begin{equation*}
  \PC_n^{\langle m \rangle} := \bigvee_{i=0}^{m-1}f_n^{-i}\PC_{i+n},\quad n\in\Z_+.%
\end{equation*}
\end{enumerate}
\end{definition}

Given an admissible class $\EC$, we can define the measure-theoretic entropy of $f_{\infty}$ w.r.t.~this class as%
\begin{equation*}
  h_{\EC}(f_{\infty}) = h_{\EC}(f_{\infty};\mu_{\infty}) := \sup_{\PC_{\infty}\in\EC}h_{\mu_{\infty}}(f_{\infty};\PC_{\infty}).%
\end{equation*}

Some elementary properties of admissible classes and their entropy are summarized in the following proposition, cf.~\cite{Ka1}.%

\begin{proposition}
Given a measure-theoretic NDS, the following statements hold:%
\begin{enumerate}
\item[(i)] There exists a maximal admissible class $\EC_{\max}$ defined as the family of all sequences $\PC_{\infty}$ satisfying Axiom (A).%
\item[(ii)] Unions and nonempty intersections of admissible classes are admissible classes.%
\item[(iii)] For each $\emptyset \neq \FC \subset \EC_{\max}$ there exists a smallest admissible class $\EC(\FC)$ containing $\FC$, and its entropy satisfies%
\begin{equation*}
  h_{\EC(\FC)}(f_{\infty}) = \sup_{\PC_{\infty}\in\FC}h(f_{\infty};\PC_{\infty}).%
\end{equation*}
\end{enumerate}
\end{proposition}

One might be tempted to regard the maximal admissible class $\EC_{\max}$ as a canonical admissible class for the definition of entropy. However, this class is usually useless, because it contains two many elements. In \cite[Ex.~18]{Ka1} it has been shown that $h_{\EC_{\max}}(f_{\infty}) = \infty$ whenever the maps $f_n$ are bi-measurable and the probability spaces $X_n$ are non-atomic.%

As in the classical theory, we can describe the dependence of $h(f_{\infty};\PC_{\infty})$ on $\PC_{\infty} \in \EC_{\max}$, using a metric on $\EC_{\max}$, defined as%
\begin{equation*}
  D(\PC_{\infty},\QC_{\infty}) := \sup_{n\in\Z_+}\left(H_{\mu_n}(\PC_n|\QC_n) + H_{\mu_n}(\QC_n|\PC_n)\right),%
\end{equation*}
with the conditional entropy $H(\cdot|\cdot)$. In the classical case, $D(\cdot,\cdot)$ reduces to the well-known \emph{Rokhlin metric}. Just as in this case, the map $\PC_{\infty} \mapsto h(f_{\infty};\PC_{\infty})$ is Lipschitz continuous w.r.t.~$D$ with Lipschitz constant $1$.%

One particularly useful property of the measure-theoretic entropy w.r.t.~an admissible class is the following power rule, cf.~\cite[Prop.~25]{Ka1}.%

\begin{proposition}
Given a measure-theoretic NDS $(X_{\infty},f_{\infty})$ and $m\in\N$, consider the $m$th power system $(X^{[m]}_{\infty},f^{[m]}_{\infty})$. If $\EC$ is an admissible class for $(X_{\infty},f_{\infty})$, we denote by $\EC^{[m]}$ the class of all sequences of partitions for $X^{[m]}_{\infty}$ which are defined by restricting the sequences in $\EC$ to the spaces in $X^{[m]}_{\infty}$, i.e., $\PC_{\infty} = \{\PC_n\}_{n\in\Z_+}\in\EC$ iff%
\begin{equation*}
  \PC^{[m]}_{\infty} := \{\PC_{nm}\}_{n\in\Z_+} \in \EC^{[m]}.%
\end{equation*}
Then $\EC^{[m]}$ is an admissible class for $(X^{[m]}_{\infty},f^{[m]}_{\infty})$ and%
\begin{equation*}
  h_{\EC^{[m]}}\left(f^{[m]}_{\infty}\right) = m\cdot h_{\EC}\left(f_{\infty}\right).%
\end{equation*}
\end{proposition}

\subsection{Measure-theoretic entropy for topological NDS}\label{subsec_mestopent}%

The concept of measure-theoretic entropy described in the preceding subsection appears to be too general and abstract for interesting applications. In this section, we explain how measure-theoretic and topological entropy interact through the definition of a specific admissible class adapted to the metric space structure of a topological NDS.%

In the following, let $(X_{\infty},f_{\infty})$ be a topological NDS and $\mu_{\infty}$ an associated IMS.%

\begin{definition}
The \emph{Misiurewicz class} $\EC_{\Mis}$ associated with $(X_{\infty},f_{\infty})$ and $\mu_{\infty}$ is defined as follows. A sequence $\PC_{\infty} = (\PC_n)_{n\in\Z_+}$ of finite Borel partitions, $\PC_n = \{P_{n,1},\ldots,P_{n,k_n}\}$, belongs to $\EC_{\Mis}$ if for every $\ep>0$ there are $\delta>0$ and compact sets $K_{n,i} \subset P_{n,i}$ for $n\in\Z_+$, $1\leq i \leq k_n$, such that the following holds for all $n\in\Z_+$:%
\begin{enumerate}
\item[(a)] $\mu_n(P_{n,i}\backslash K_{n,i}) \leq \ep$ for $1\leq i\leq k_n$.%
\item[(b)] If $x\in K_{n,i}$, $y\in K_{n,j}$, $i\neq j$, then $d_n(x,y) \geq \delta$.%
\end{enumerate}
\end{definition}

As it turns out, this definition in fact yields an admissible class that is well-adapted to the metric space structure, as expressed by the following theorem.%

\begin{theorem}
$\EC_{\Mis}$ is an admissible class with the following properties:%
\begin{enumerate}
\item[(i)] $\EC_{\Mis}$ and the associated entropy $h_{\EC_{\Mis}}(f_{\infty};\mu_{\infty})$ are preserved by equi-conjugacies, i.e., equicontinuous changes of coordinates.%
\item[(ii)] In the autonomous case, i.e., when $X_{\infty},d_{\infty},f_{\infty}$ and $\mu_{\infty}$ are constant, $h_{\EC_{\Mis}}(f_{\infty};\mu_{\infty})$ reduces to the usual Kolmogorov-Sinai measure-theoretic entropy.%
\item[(iii)] The inequality%
\begin{equation*}
  h_{\EC_{\Mis}}(f_{\infty};\mu_{\infty}) \leq h_{\tp}(f_{\infty})%
\end{equation*}
holds (establishing one part of the variational principle).%
\end{enumerate}
\end{theorem}

The proofs of (i) and (iii) can be found in \cite[Prop.~26, Prop.~27, Thm.~28]{Ka1} and the proof of (ii) in \cite[Cor.~3.1]{KLa}.%

Since the definition of $\EC_{\Mis}$ is tailored to the (first half of the) proof of the variational principle due to Misiurewicz \cite{Mis}, proving (ii) is an easy task. However, it is not as easy as it might seem to prove that $h_{\EC_{\Mis}}$ in fact generalizes the classical notion of measure-theoretic entropy, since even if $X_{\infty}$, $d_{\infty}$, $f_{\infty}$ and $\mu_{\infty}$ are assumed to be constant, we still have to deal with non-constant sequences of partitions. The proof is accomplished through the following result, cf.~\cite[Thm.~3.1]{KLa}.%

\begin{theorem}\label{thm_finescaleseq}
Assume that there exists a sequence $(\RC_{\infty}^k)_{k\in\Z_+}$ in $\EC_{\Mis}$ with%
\begin{equation*}
  \lim_{k\rightarrow\infty}\sup_{n\in\Z_+}\sup_{R \in \RC^k_n}\diam R = 0.%
\end{equation*}
Then the measure-theoretic entropy satisfies%
\begin{equation*}
  h_{\EC_{\Mis}}(f_{\infty};\mu_{\infty}) = \lim_{k\rightarrow\infty}h(f_{\infty};\RC^k_{\infty}) = \sup_{k\in\Z_+}h(f_{\infty};\RC^k_{\infty}).%
\end{equation*}
\end{theorem}

In the autonomous case, it is clear that every constant sequence of partitions is contained in $\EC_{\Mis}$, hence any refining sequence of partitions defines a sequence $(\RC_{\infty}^k)_{k\in\Z_+}$, as required in the theorem. Consequently, the theorem says that the entropy is already determined on the constant sequences of partitions, so the classical definition of Kolmogorov-Sinai entropy is retained.%

In general, it is unclear whether the Misiurewicz class contains sequences as required in Theorem \ref{thm_finescaleseq}. The following result, proved in \cite{KLa}, yields several sufficient conditions in the case when the state space is time-invariant, cf.~\cite[Thm.~3.2]{KLa}.%

\begin{theorem}\label{thm_finescaleconds}
Assume that $(X_n,d_n) \equiv (X,d)$ for some compact metric space $(X,d)$. Then each of the following conditions guarantees that $\EC_{\Mis}$ contains elements of arbitrarily (uniformly) small diameter:%
\begin{enumerate}
\item[(i)] $\{\mu_n : n\in\Z_+\}$ is relatively compact in the strong topology on the space of measures.%
\item[(ii)] For every $\alpha>0$ there is a finite measurable partition $\AC$ of $X$ with $\diam\AC<\alpha$ such that $\nu(\partial\AC)=0$ for all weak$^*$-limits $\nu$ of $\mu_{\infty}$. (This holds, in particular, if there are only countably many non-equivalent weak$^*$-limits.)%
\item[(iii)] $X = [0,1]$ or $X = \rmS^1$ and there exists a dense set $D\subset X$ such that every $x\in D$ satisfies $\nu(\{x\})=0$ for all weak$^*$-limits $\nu$ of $\mu_{\infty}$.%
\item[(iv)] $X$ has topological dimension zero.%
\end{enumerate}
In each case, the sequences of partitions can in fact be chosen constant.%
\end{theorem}

The following theorem provides an example, where both topological and measure-theoretic entropy can be computed, cf.~\cite[Thm.~5.4 and Thm.~5.5]{Ka2}.%

\begin{theorem}
Let $M$ be a compact Riemannian manifold and $f_{\infty} = (f_n)_{n\in\Z_+}$ a sequence of $C^2$-expanding maps $f_n:M \rightarrow M$ with expansion factors uniformly bounded away from one, and $C^2$-norms uniformly bounded. Then%
\begin{equation*}
  h_{\tp}(f_{\infty}) = \limsup_{n\rightarrow\infty}\frac{1}{n}\log \int_M |\det\rmD f_0^n(x)| \rmd \vol,%
\end{equation*}
and for any smooth initial measure $\mu_0$, with $\mu_{\infty} = (f_0^n\mu_0)_{n\in\Z_+}$,%
\begin{equation*}
  h_{\EC_{\Mis}}(f_{\infty};\mu_{\infty}) = \limsup_{n\rightarrow\infty}\frac{1}{n}\int_M \log|\det\rmD f_0^n(x)| \rmd \vol.%
\end{equation*}
\end{theorem}

The question under which conditions an NDS satisfies a full variational principle, i.e.,%
\begin{equation*}
  h_{\tp}(f_{\infty}) = \sup_{\mu_{\infty}}h_{\EC_{\Mis}}(f_{\infty};\mu_{\infty})%
\end{equation*}
is completely open. Only some examples are known which do not allow for a broad generalization.%

\section{An example}\label{sec_ex}%

In this section, we apply the theory explained above to an NDS which has been introduced in \cite{BOp} by Balibrea and Oprocha. We will need the following proposition whose proof is completely analogous to the autonomous case, and hence is omitted.%

\begin{figure}[b]
\begin{center}
  \subfigure[Graph of $f$]{\includegraphics[scale=.65]{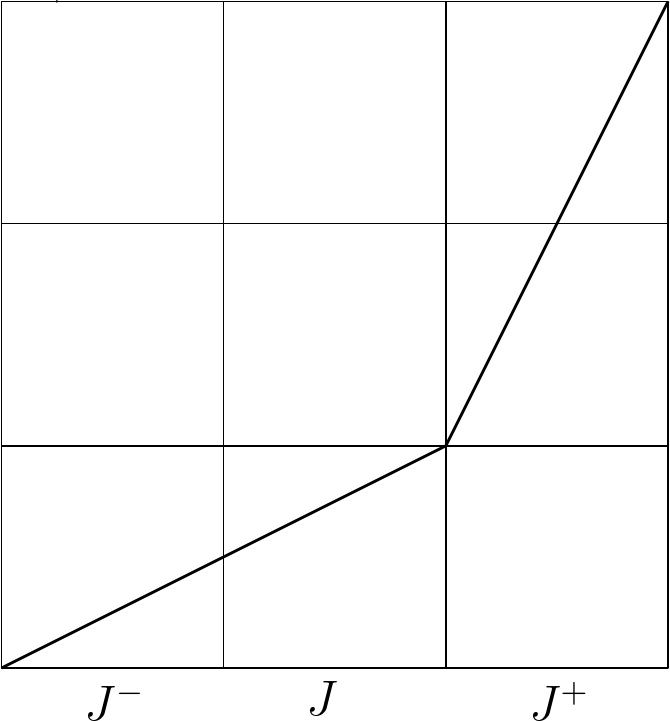}}
  \subfigure[Graph of $g$]{\includegraphics[scale=.65]{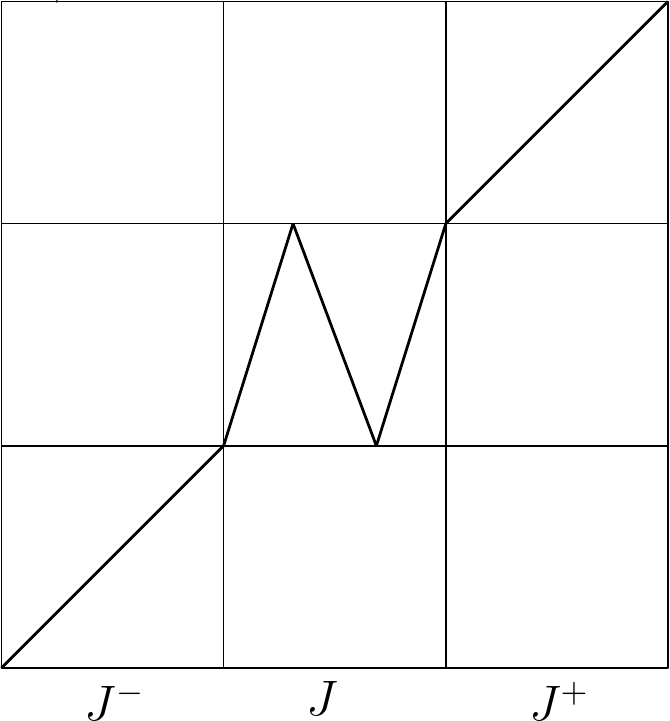}}
\end{center}
\caption{The maps $f$ and $g$}
\label{fig:1}
\end{figure}

\begin{proposition}\label{prop_lipest}
Let $(X_{\infty},f_{\infty})$ be a topological NDS such that $f_n$ is (globally) Lipschitz-continuous with Lipschitz-constant $L_n$ for each $n$ and $X_0$ has finite upper capacitive dimension $\overline{\dim}_C(X_0)$. Then%
\begin{equation*}
  h_{\tp}(f_{\infty}) \leq \overline{\dim}_C(X_0) \cdot \limsup_{n\rightarrow\infty}\frac{1}{n}\sum_{i=0}^{n-1}\max\{0,\log L_i\}.%
\end{equation*}
\end{proposition}

Now consider the NDS from \cite[Thm.~4]{BOp}, which is constructed from the two piecewise affine maps depicted in Fig.~\ref{fig:1}. More precisely, let $m_0 := 1$ and $m_n := 2^{n^2}$ for all $n\in\N$. Consider the maps $f,g:[0,1]\rightarrow[0,1]$ in Fig.~\ref{fig:1}, and the NDS $f_{\infty} = (f_n)_{n\in\Z_+}$ defined by%
\begin{equation*}
  f_i := \left\{\begin{array}{rl}
	          f & \mbox{ if } i = m_n \mbox{ for some } n\\
						g & \mbox{ otherwise}
					\end{array}\right..%
\end{equation*}
For the Lebesgue measure $\lambda$ on $[0,1]$ we have weak convergence $\mu_n = f_0^n\lambda \rightarrow \delta_0$, since every trajectory with initial value in $[0,1)$ converges to zero. More precisely, this implies $\varphi \circ f_0^n(x) \rightarrow \varphi(0)$ for every $x\in [0,1)$ and every continuous function $\varphi:[0,1]\rightarrow\R$. Hence, $\int \varphi \rmd\mu_n = \int \varphi \circ f_0^n \rmd\lambda \rightarrow \int \varphi(0) \rmd \lambda$ by the theorem of dominated convergence. Consequently, by Theorem \ref{thm_finescaleconds}(ii), the admissible class $\EC_M(\mu_{\infty})$ contains all constant sequences of partitions with $\delta_0$-zero boundaries, in particular all constant sequences $\PC_n \equiv \PC$, where $\PC$ consists of nontrivial subintervals of $[0,1]$.%

Let $\PC$ be a partition of $[0,1]$ into intervals of length $1/(3k)$ for some $k\in\N$. Then each interval in $\PC$ is completely contained in $J^- := [0,1/3]$, $J := [1/3,2/3]$ or $J^+ := [2/3,1]$. Let $\lambda$ denote the Lebesgue measure on $[0,1]$. Then%
\begin{equation*}
  H_{\lambda}\left(\bigvee_{i=0}^{m_n}f_0^{-i}\PC\right) = H_{\lambda}\left(\bigvee_{i=0}^{m_{n-1}}f_0^{-i}\PC \vee \bigvee_{i=m_{n-1} + 1}^{m_n}f_0^{-i}\PC\right)
	\geq H_{\lambda}\left(\bigvee_{i=m_{n-1}+1}^{m_n}f_0^{-i}\PC\right).%
\end{equation*}
Note that for $m_{n-1}+1 \leq i \leq m_n$ we have%
\begin{equation*}
  f_0^{-i} = \left(g^{i - m_{n-1} - 1} \circ f_{m_{n-1}} \circ \cdots \circ f_1 \circ f_0 \right)^{-1} = f_0^{-(m_{n-1}+1)} \circ g^{-(i-m_{n-1}-1)},%
\end{equation*}
and hence, writing $l_n := m_n - m_{n-1} - 1$,%
\begin{equation*}
  H_{\lambda}\left(\bigvee_{i=0}^{m_n}f_0^{-i}\PC\right) \geq H_{\lambda}\left(f_0^{-(m_{n-1}+1)}\bigvee_{i=0}^{l_n} g^{-i}\PC\right).%
\end{equation*}
Now we look only at those members of $\bigvee_{i=0}^{l_n}g^{-i}\PC$ that come from intervals $P \in \PC$ with $P \subset J$. Let us write $\PC^J$ for the the set of all elements in $\PC$ contained in $J$. Then the above can be estimated by%
\begin{align*}
  &\geq H_{\lambda}\left(f_0^{-(m_{n-1}+1)}\bigvee_{i=0}^{l_n} g^{-i} \PC^J \right)\\
	&= - \sum_{P \in \bigvee_{i=0}^{l_n}g^{-i}\PC^J} \lambda(f_0^{-(m_{n-1}+1)}P) \log \lambda(f_0^{-(m_{n-1}+1)}P).%
\end{align*}
Now we use that $J$ is $g$-invariant and $f_0^{-(m_{n-1}+1)}(A) = f^{-n}(A)$ for any $A \subset J$ and $n\geq1$. Moreover, we use that $f^{-1}(x) = (1/2)(x - (1/3)) + (2/3)$ on $J$. Together with the fact that $g^{-1}$ is trivial on $J^+$, this gives%
\begin{align*}
  H_{\lambda}\left(\bigvee_{i=0}^{m_n}f_0^{-i}\PC\right) &\geq - \left(\# \bigvee_{i=0}^{l_n}g^{-i}\PC^J\right) \frac{1}{3^{l_n}2^n 3k} \log \frac{1}{3^{l_n}2^n 3k}\allowdisplaybreaks\\
	&= \log \left(3^{l_n}2^n 3k\right) = l_n\log (3) + n\log(2) + \log(3k).%
\end{align*}
Dividing by $m_n$ and sending $n$ to infinity, gives $\log(3)$, since%
\begin{equation*}
  \frac{m_n - m_{n-1} - 1}{m_n} = 1 - 2^{- 2n - 1} - \frac{1}{2^{n^2}} \rightarrow 1,%
\end{equation*}
and $n/m_n \rightarrow 0$. Writing $\lambda_{\infty}$ for the sequence $\lambda_n := f_0^n\lambda$, we obtain%
\begin{equation*}
  h_{\EC_M}(f_{\infty};\lambda_{\infty}) \geq \log(3).%
\end{equation*}
Since $L=3$ is a Lipschitz constant for both $f$ and $g$, Proposition \ref{prop_lipest} yields%
\begin{equation*}
  \log(3) \leq h_{\EC_M}(f_{\infty};\lambda_{\infty}) \leq h_{\tp}(f_{\infty}) \leq \log(3),%
\end{equation*}
implying that for $f_{\infty}$ a full variational principle is satisfied with $\lambda_{\infty}$ being an IMS of maximal entropy.%

\begin{remark}
It is easy to see that every trajectory $\{f_0^n(x)\}_{n\in\Z_+}$ with $x \neq 1$ converges to $0$. Hence, the example shows that both the measure-theoretic and the topological entropy can capture transient chaotic behavior, which is not seen in the asymptotic behavior of trajectories.%
\end{remark}


\begin{thebibliography}{999}
\bibitem{AKM} R.~L.~Adler, A.~G.~Konheim, M.~H.~McAndrew. \emph{Topological entropy}. Trans. Am. Math. Soc. 114 (1965), 309--319.%
\bibitem{BOp} F.~Balibrea, P.~Oprocha. \emph{Weak mixing and chaos in nonautonomous discrete systems}. Appl. Math. Lett. 25 (2012), no. 8, 1135--1141.%
\bibitem{Bow} R.~Bowen. \emph{Entropy for group endomorphisms and homogeneous spaces}. Trans. Am. Math. Soc. 153 (1971), 401--414.%
\bibitem{Can} J.~S.~C\'anovas. \emph{On entropy of nonautonomous discrete systems}. Progress and Challenges in Dynamical Systems. Springer (2013), 143--159.%
\bibitem{Din} E.~I.~Dinaburg. \emph{The relation between topological entropy and metric entropy}. Dokl.\ Akad.\ Nauk SSSR 190, 19--22 (1970), (Soviet Math.\ Dokl.\ 11 (1969), 13--16).%
\bibitem{Dow} T. Downarowicz. Entropy in Dynamical Systems. New Mathematical Monographs 18. Cambridge University Press, Cambridge, 2011.%
\bibitem{Go1} T.~N.~T.~Goodman. \emph{Relating topological entropy and measure entropy}. Bull.\ Lond.\ Math.\ Soc.\ 3 (1971), 176--180.%
\bibitem{Go2} T.~N.~T.~Goodman. \emph{Topological sequence entropy}. Proc. London Math. Soc. (3) 29 (1974), 331--350.%
\bibitem{Gow} L.~W.~Goodwyn. \emph{Topological entropy bounds measure-theoretic entropy}. Proc.\ Am.\ Math.\ Soc.\ 23 (1969), 679--688.%
\bibitem{Ka1} C.~Kawan. \emph{Metric entropy of nonautonomous dynamical systems}. Nonauton.\ Stoch.\ Dyn.\ Syst.\ 1 (2013), 26--52.%
\bibitem{Ka2} C.~Kawan. \emph{Expanding and expansive time-dependent dynamics}. Nonlinearity 28 (2015), no. 3, 669--695.%
\bibitem{Ka3} C.~Kawan. \emph{Exponential state estimation, entropy and Lyapunov exponents}. Submitted, 2016; arXiv:1605.03210 [math.DS]%
\bibitem{KLa} C.~Kawan, Y.~Latushkin. \emph{Some results on the entropy of non-autonomous dynamical systems}. Dyn. Syst. 31 (2016), no. 3, 251--279.%
\bibitem{KSn} S.~Kolyada, L.~Snoha. \emph{Topological entropy of nonautonomous dynamical systems}. Random Comput. Dynamics 4 (1996), no. 2--3, 205--233.%
\bibitem{KMS} S.~Kolyada, M.~Misiurewicz, L.~Snoha. \emph{Topological entropy of nonautonomous piecewise monotone dynamical systems on the interval}. Fund. Math. 160  (1999), no. 2, 161--181.%
\bibitem{MPo} A.~S.~Matveev, A.~Pogromsky. \emph{Observation of nonlinear systems via finite capacity channels: Constructive data rate limits}. Automatica 70 (2016), 217--229.%
\bibitem{Orn} D.~S.~Ornstein. \emph{Bernoulli shifts with the same entropy are isomorphic}. Advances in Mathematics, 4 (1970), 337--352.%
\bibitem{PMa} A.~Y.~Pogromsky, A.~S.~Matveev. \emph{Estimation of topological entropy via the direct Lyapunov method}. Nonlinearity 24 (2011), no. 7, 1937--1959.%
\bibitem{Mis} M.~Misiurewicz. \emph{Topological entropy and metric entropy}. Ergodic theory (Sem., Les Plans-sur-Bex, 1980) (French), 61--66, Monograph. Enseign. Math., 29, Univ. Gen\'eve, Geneva (1981).%
\bibitem{Mou} C.~Mouron. \emph{Positive entropy on nonautonomous interval maps and the topology of the inverse limit space}. Topology Appl. 154 (2007), no. 4, 894--907.%
\bibitem{OW1} P.~Oprocha, P.~Wilczynski. \emph{Chaos in nonautonomous dynamical systems}. An. Stiint. Univ. ``Ovidius'' Constanta Ser. Mat. 17 (2009), no. 3, 209--221.%
\bibitem{OW2} P.~Oprocha, P.~Wilczynski. \emph{Topological entropy for local processes}. J. Differential Equations 249 (2010), no. 8, 1929--1967.%
\bibitem{ZCh} J.~Zhang, L.~Chen. \emph{Lower bounds of the topological entropy for nonautonomous dynamical systems}. Appl. Math. J. Chinese Univ. Ser. B 24 (2009), no. 1, 76--82.%
\bibitem{ZZH} Y.~Zhu, J.~Zhang, L.~He. \emph{Topological entropy of a sequence of monotone maps on circles}. Korean Math. Soc. 43 (2006), no. 2, 373--382.%
\bibitem{ZLX} Y.~Zhu, Z.~Liu, X.~Xu, W.~Zhang. \emph{Entropy of nonautonomous dynamical systems}. J. Korean Math. Soc. 49 (2012), no. 1, 165--185.%
\end{thebibliography}
\end{document}